\documentclass[11pt]{article}

\usepackage{amsmath}
\usepackage{amssymb}
\usepackage[all]{xy}
\usepackage{mathrsfs}

\DeclareMathAlphabet{\scr}{U}{eus}{m}{n}

\newcommand\Z{{\mathbb Z}}
\newcommand\Q{{\mathbb Q}}
\newcommand\F{{\mathbb F}}

\newcommand\cB{{\cal B}}

\newcommand\cF{{\cal F}}

\newcommand\Tr{{\rm Tr}}
\newcommand\Frob{{\rm Fr}}
\newcommand\Fr{{\rm Fr}}
\newcommand\xt{[\hspace{-.12em}[ t ]\hspace{-.12em} ] } 
\newcommand\xT{(\!( t )\!)}

\newcommand\Ql{\overline{\mathbb Q}_\ell}

\newcommand\Fl{\mathcal Fl}

\newcommand\tens\otimes

\newcommand{\Grass}{\mathop{\rm Grass}\nolimits}

\newcommand\qed{\hfill$\square$}

\newtheorem{thm}{Theorem}[section]
\newtheorem{stz}[thm]{Proposition}
\newtheorem{lem}[thm]{Lemma}
\newtheorem{Def}[thm]{Definition}
\newtheorem{kor}[thm]{Corollary}

\numberwithin{equation}{section}

\title{Bounds on weights of nearby cycles and Wakimoto sheaves on affine
flag manifolds} \author{Ulrich G\"{o}rtz and Thomas J. Haines
\footnote{Research of the second author partially supported by NSF grant
DMS 0303605 and a Sloan Research Fellowship.}} \date{}

\begin{document}

\maketitle

\section{Introduction}

Let $G$ be a split connected reductive group over a finite field $\F_p$
with algebraic closure $k$, fix an Iwahori subgroup $\cB \subset
G({\F}_p(\xt))$ and let $\Fl = G(k \xT )/\cB_k$ denote the affine flag
variety of $G$.  Let $q$ denote a power of $p$ and fix a prime $\ell \neq
p$.  In \cite{GH}, the authors study the Jordan-H\"{o}lder series for
objects in the Hecke category $P^\cB_q(\Fl, \Ql)$.  This is the category of
$\cB$-equivariant perverse Weil-sheaves $\cF$ on $\Fl$ having the property
that for any $x \in \Fl(\F_q)$, the trace of the Frobenius $\Fr_q$ on the
stalk $\cF_x$ satisfies $$ \Tr(\Fr_q, \cF_x) \in \Z[q,q^{-1}].  $$ Put more
precisely, for each $x$ there exists a polynomial $P(X,Y) \in \Z[X,Y]$ such
$\Tr(\Fr_{q^n}, \cF_x) = P(q^{n},q^{-n})$, for all positive integers $n$.  

By \cite{GH}, $\S 4$ the objects in $P^\cB_q(\Fl,\Ql)$ are mixed and have
finite length, and the irreducible objects are the integral Tate-twists
$IC_w(i)$ of the intersection complexes $IC_w = j_{w,!*}(\Ql)[\ell(w)]$.
Here for $w$ belonging to the extended affine Weyl group $\widetilde{W}$ we
take the Goresky-MacPherson middle extension along the immersion $j_w:
\cB_w \hookrightarrow \cB$ of the corresponding Schubert cell $\cB_w$.  We
may thus define non-negative multiplicities $m(\cF,w,i)$ by the identity in
$P^{\cB}_q(\Fl,\Ql)$ $$ \cF^{ss} = \bigoplus_{w \in \widetilde{W}}
\bigoplus_{i \in \Z} IC_w (-i)^{\oplus m(\cF, w, i)}.  $$ We formally
encode the data of the integers $m(\cF,w,i)$ with the {\em multiplicity
function} $$ m(\cF,w) := \sum_i m(\cF,w,i) q^i \in \Z[q,q^{-1}].  $$

In \cite{GH} we studied the functions $m(R\Psi,w)$ for certain nearby
cycles sheaves $R\Psi$ on $\Fl$ which arise as follows.  If $G$ is the
group ${\rm GL}_n$ or ${\rm GSp}_{2n}$, there is a $\Z_p$-ind-scheme $M$
which is a deformation of the affine Grassmannian $\Grass_{\mathbb Q_p}$ to
the affine flag variety $\Fl_{\F_p}$ for the underlying group $G$.  The
ind-scheme $M$ has been constructed over $\Z_p$ for these two groups in
\cite{HN1}, and it is related to the bad reduction of certain Shimura
varieties (see below).  In the function field case the analogous
deformation ${\rm FL}_X$ over a curve $X$ exists for all split connected
reductive groups, and has been studied by Gaitsgory \cite{Ga}.  In that
case, there is a distinguished point $x_0 \in X$, and the fiber of ${\rm
FL}_X$ over $x \in X$ is ${\Fl}_k$ if $x = x_0$ and is ${\rm Grass}_k
\times G/B$ if $x \neq x_0$, where $G/B$ is the flag variety which is the
``reduction modulo $t$'' of $G(k\xt)/\cB_k$.  Let ${\mathcal Q}_\mu \subset
\Grass$ denote the stratum indexed by a dominant coweight $\mu$ of $G$, and
let $IC_\mu$ denote the intersection complex $j_{\mu,!*}(\Ql)[{\rm dim}\,
{\mathcal Q}_\mu]$.  Let $R\Psi_\mu = R\Psi(IC_\mu)$, the nearby cycles
being taken with respect to the deformation $M$ or ${\rm FL}_X$.  More
precisely, in the $p$-adic setting we write $R\Psi_\mu = R\Psi^M(IC_\mu)$,
and in the function-field setting we write $R\Psi_\mu = R\Psi^{{\rm
FL}_X}(IC_\mu \boxtimes \delta)$, where $\delta$ is the skyscraper sheaf
supported at the base point of $G/B$.  

The sheaf $R\Psi_\mu$ belongs to $P^{\cB}_q(\Fl,\Ql)$ (see \cite{GH},$\S
2,6$), a fact that is essentially a consequence of the Kottwitz conjecture,
see (\ref{eq:Kottwitz_conj}).  Hence we may consider the multiplicity
functions $m(R\Psi_\mu,w)$ for $w \in \widetilde{W}$.  By \cite{GH}, $\S 7$
it is known that $m(R\Psi_\mu,w) \neq 0$ if and only $w \in {\rm
Adm}(\mu)$, the finite subset of $\mu$-admissible elements in
$\widetilde{W}$ (see \cite{KR},\cite{HN2} for explicit descriptions of
${\rm Adm}(\mu)$).  

Our main theorem is the following result conjectured in \cite{GH} and
proved there for the case $G = {\rm GL}_n$ or $\mu$ minuscule. 

\begin{thm} \label{nearby.poly}  For any $w \in {\rm Adm}(\mu)$, the
function $m(R\Psi_\mu,w)$ is a polynomial in $q$ having degree at most
$\ell(\mu) - \ell(w)$.  Equivalently, the same statement holds for the
functions $\Tr(\Fr_q,(R\Psi_{\mu})_w)$.  Equivalently, for each closed
point $z \in \cB_w$, the weights of the complex
$R\Psi_{\mu}(\frac{\ell(\mu)}{2})_z$ are $\leq \ell(\mu)-\ell(w)$.
\end{thm} Here the length function $\ell: \widetilde{W} \rightarrow
\Z_{\geq 0}$ is that defined by the simple affine reflections $S_{\rm aff}$
through the walls of the alcove in the Bruhat-Tits building which is fixed
by the Iwahori $\cB$ (see \cite{GH}).  The symbol $\ell(\mu)$ denotes the
length of the translation element $t_\mu \in \widetilde{W}$ corresponding
to the cocharacter $\mu$.

The study of $m(R\Psi,w)$ relates to nearby cycles on certain  Shimura
varieties with Iwahori level structure at $p$.  In fact, for such a Shimura
variety $Sh = Sh({\bf G},h)_{\bf K}$ where the $p$-adic group ${\bf
G}_{\Q_p}$ is either ${\rm GL}_n \times {\mathbb G}_m$ or ${\rm GSp}_{2n}$,
the Rapoport-Zink local model ${\bf M}^{\rm loc}$ for the singularities in
the special fiber $Sh_{\overline{\F}_p}$ can be realized as the
scheme-theoretic closure in $M$ of ${\mathcal Q}_\mu \subset {\rm
Grass}_{\Q_p}$, where $\mu$ is the minuscule coweight coming from the
Shimura datum $h$.  See \cite{GH}, $\S 8$ for a more detailed discussion of
the connection with Shimura varieties.

\medskip

A related result concerns the {\em Wakimoto sheaves} $M_{u,v}$ on $\Fl$.
Here for elements $u,v \in \widetilde{W}$ we define following I. Mirkovi\'c
the sheaf $$ M_{u,v} = j_{u!}(\Ql)[\ell(u)] * j_{v*}(\Ql)[\ell(v)], $$
using the standard convolution operation $*: P^{\cB}_q(\Fl,\Ql) \times
P^{\cB}_q(\Fl,\Ql) \rightarrow D^{b,\rm Weil}_c(\Fl,\Ql)$, the latter being
the category consisting of objects of the ``derived category of $\ell$-adic
sheaves'' $D^b_c(\Fl,\Ql)$, equipped with a Weil structure.  A result of
Mirkovi\'c states that $M_{u,v}$ is perverse, hence it belongs to
$P^{\cB}_q(\Fl,\Ql)$ (Mirkovi\'c's proof appears in \cite{HP} and in
\cite{AB}).  Let us renormalize and define $$ \widetilde{M}_{u,v} =
M_{u,v}(\frac{\ell(u) + \ell(v) - \ell(uv)}{2}).$$  It is known that the
function $\Tr(\Fr_q, M_{u,v})$ is supported on a subset of the elements $\{
x \in \widetilde{W} ~ | ~ x \leq uv \}$, where the Bruhat order $\leq$ is
defined using $S_{\rm aff}$ (see \cite{H01}).

\begin{thm} \label{wakimoto.poly} For any $w \leq uv$, the function
$m(\widetilde{M}_{u,v},w)$ is a polynomial in $q$ having degree at most
$\ell(uv) - \ell(w)$.  Equivalently, the same holds for the functions
$\Tr(\Fr_q, (\widetilde{M}_{u,v})_w)$.  Equivalently, for each closed point
$z \in \cB_w$, the weights of the complex
$\widetilde{M}_{u,v}(\frac{\ell(uv)}{2})_z$ are $\leq \ell(uv) - \ell(w)$.
\end{thm}

As we explain at the end of section \ref{wakimoto_section}, Theorem
\ref{wakimoto.poly} can be used to prove Theorem \ref{nearby.poly}.  But in
section \ref{nearby_section} we will give a more direct proof of Theorem
\ref{nearby.poly} using an approach we introduced (but did not complete) in
\cite{GH}, $\S 10$. 

Theorem \ref{wakimoto.poly} might be interesting in its own right, since
Wakimoto sheaves have played an important role in several recent papers
which study the geometry relating the two presentations (Iwahori-Matsumoto
and Bernstein) of an affine Hecke algebra.  See \cite{HP}, and especially
\cite{AB} and \cite{ABG}.

\section{From bounds on weights to bounds on degrees}
 
Let ${\rm supp}(\cF) = \{ x \in \widetilde{W} ~ | ~ \Tr(\Fr_q,\cF_x) \neq 0
\}$, and let us define the $IC$-{\em support} to be the finite set ${\rm
ICsupp}(\cF)= \{ x \in \widetilde{W} ~ | ~ m(\cF, x) \neq 0 \}$.  An
argument in \cite{GH}, $\S7$ shows that ${\rm supp}(\widetilde{M}_{u,v})
\subseteq \{ x \in \widetilde{W} ~ | ~ x \leq uv \} = {\rm
ICsupp}(\widetilde{M}_{u,v})$.  Using the Kottwitz conjecture (see end of
section \ref{wakimoto_section}), we get the analogous result for nearby
cycles: ${\rm supp}(R\Psi_\mu) \subseteq {\rm Adm}(\mu) = {\rm
ICsupp}(R\Psi_\mu)$.  For general $\cF$, the precise relation between the
sets ${\rm supp}(\cF)$ and ${\rm ICsupp}(\cF)$ is not clear.  However, from
equation (\ref{eq:recursion_formula}) below, it is clear that any maximal
element of ${\rm ICsupp}(\cF)$ belongs to ${\rm supp}(\cF)$, and
vice-versa.  Hence the two sets have the same maximal elements.

Let us recall a definition from \cite{GH}.  Fix an integer $d$.  We say
$\cF$ satisfies property $(P)_d$ if for all $y \in \widetilde{W}$ we have
the equivalent identities \begin{align} \label{eq:propertyP} \Tr(\Fr_q,
\cF_y) &= \varepsilon_d \varepsilon_y q^{d} q_y^{-1}
\overline{\Tr(\Fr_q,\cF_y)} \\ \Tr(\Fr_q, D\cF_y) &= \varepsilon_d
\varepsilon_y q^{-d} q_y^{-1} \overline{\Tr(\Fr_q,D\cF_y)}, \end{align}
where $\varepsilon_d = (-1)^d$, $\varepsilon_y = (-1)^{\ell(y)}$, $q_y =
q^{\ell(y)}$, $D\cF$ denotes the Verdier dual of $\cF$, and $h \mapsto
\overline{h}$ is the involution of $\Z[q^{1/2},q^{-1/2}]$ determined by
$q^{1/2} \mapsto q^{-1/2}$.  For example, in \cite{GH} we proved that
$R\Psi_\mu$ satisfies $(P)_d$ with $d = \ell(\mu)$ and
$\widetilde{M}_{u,v}$ satisfies $(P)_d$ with $d = \ell(uv)$.   

\begin{lem} \label{weights_v_degrees}  Let $\cF$ be any object of
$P^{\cB}_q(\Fl,\Ql)$ and let $d$ be any positive integer.  Then for any $x
\in \widetilde{W}$, the following statements are equivalent:
\begin{enumerate} \item [(1)] For all $w \geq x$, $\Tr(\Fr_q,\cF_w)$ is a
polynomial in $q$ having degree $\leq d - \ell(w)$; \item [(2)] For all $w
\geq x$, $m(\cF,w)$ is a polynomial in $q$ having degree $\leq d -
\ell(w)$.  \end{enumerate} If $\cF$ satisfies property $(P)_d$, then these
are also equivalent to 

~~(3) For all $w \geq x$, the weights of the stalk complex
$\cF(\frac{d}{2})_w$ are $\leq d - \ell(w)$.  \end{lem}

\noindent{\em Proof.} For $x,y \in \widetilde{W}$ let $P_{x,y}(q)$ denote
the corresponding Kazhdan-Lusztig polynomial.  For any $x \in
\widetilde{W}$ we have using \cite{KL2} the identity \begin{equation}
\label{eq:recursion_formula} \Tr(\Fr_q, \cF_x) = \varepsilon_x  m(\cF,x) +
\sum_{\underset{w > x}{w \in {\rm ICsupp}(\cF)}} \varepsilon_w  m(\cF,w)
P_{x,w}(q).  \end{equation} The maximal elements of the sets ${\rm
supp}(\cF)$ and ${\rm ICsupp}(\cF)$ agree, and (\ref{eq:recursion_formula})
implies the equivalence $(1) \Leftrightarrow (2)$ holds for $x$ a maximal
element.  The case of general $x$ follows by descending induction with
respect to $\,\,\leq$, using (\ref{eq:recursion_formula}) and the bound
${\rm deg}_q P_{x,w} < \frac{1}{2}(\ell(w) - \ell(x))$ if $x < w$.

Next we will prove that a statement closely related to the implications
$(3) \Leftrightarrow (2)$ holds without assuming $\cF$ satisifies property
$(P)_d$.  We claim that the following statements are equivalent.
\begin{enumerate} \item [(i)] For all $w \geq x$, the weights of
$\cF(\frac{d}{2})_w$ are $\leq d - \ell(w)$; \item[(ii)] for all $w \geq
x$, ${\rm deg}_q m(\cF,w) \leq d - \ell(w)$.  \end{enumerate}

\noindent {\em Proof of claim.}  Let us prove $(i) \Rightarrow (ii)$.  We
fix an element $w \geq x$ and we may assume $w \in {\rm ICsupp}(\cF)$.  The
condition in (i) holds for $w$ if and only if for all $w' \geq w$ and $i
\in \Z$ such that $m(\cF,w',i) \neq 0$ we have ${\rm
weight}(IC_{w'}(\frac{\ell(w')}{2})(\frac{d}{2} - \frac{\ell(w')}{2} - i))
\leq d - \ell(w)$.  Since $IC_{w'}(\frac{\ell(w')}{2})$ is a pure weight
zero perverse sheaf, this last inequality holds if and only if $2i - d +
\ell(w') \leq d - \ell(w)$, i.e., $i \leq d - \frac{\ell(w') +
\ell(w)}{2}$.  Taking $w' = w$, we get ${\rm deg}_q m(\cF,w) \leq d -
\ell(w)$, as desired.  

Now we prove $(ii) \Rightarrow (i)$.  Assume $w \geq x$.  The hypothesis
means that for every $w' \geq w$, and $i$ such that $m(\cF,w',i) \neq 0$ we
have $i \leq d - \ell(w')$ which is a fortiori $\leq d - \frac{\ell(w) +
\ell(w')}{2}$.  Now reversing the steps in the paragraph above, we see this
implies that ${\rm weight}(\cF(\frac{d}{s})_w) \leq d - \ell(w)$.  This
proves the claim.

In particular we see that (1) and (2) each imply (3) with no extra
assumptions on $\cF$.

It remains to prove that if $\cF$ satisfies $(P)_d$ and also satisfies the
equivalent conditions (i) and (ii) above, then for all $w \geq x$ the
functions $\Tr(\Fr_q,\cF_w)$ are polynomials in $q$ having degree $\leq d -
\ell(w)$.  First of all the inequality ${\rm deg}_q \Tr(\Fr_q,\cF_w) \leq d
- \ell(w)$ follows from (ii) and (\ref{eq:recursion_formula}).  The fact
that $\Tr(\Fr_q,\cF_w)$ is actually a polynomial in $q$ follows from this
inequality and the identity (\ref{eq:propertyP}) defining property $(P)_d$.
\qed

\section{A bound on the weights of nearby cycles on each stratum}
\label{nearby_section}

By Lemma \ref{weights_v_degrees}, Theorem \ref{nearby.poly} can be
reformulated as follows. 

\begin{thm} \label{nearby_bounds} For each $w \in {\rm Adm}(\mu)$ and each
closed point $z \in \cB_w$, the weights of the stalk complex
$R\Psi_{\mu}(\frac{\ell(\mu)}{2})_z$ are $\leq \ell(\mu) - \ell(w)$.
\end{thm}

To prove this we will first prove a weaker statement that holds in much
greater generality. 

% The main tools are the theories of alterations \cite{deJ} and perverse
% sheaves \cite{BBD}.

\subsection{Bounding weights generically on strata for nearby cycles over
Henselian traits}

Let $(S, s,\eta)$ be the spectrum of a complete discrete valuation ring
with finite residue field $k(s)$.  Let $k(\overline{\eta})$ be a separable
closure of the generic point $k(\eta)$, and let $\overline{S}$ be the
normalization of $S$ in $\overline{\eta}$, with closed point
$\overline{s}$.  Let $X$ be an $S$-variety, i.e. an integral separated
scheme, flat and of finite type over $S$.

Let $X_{\bar s}$ denote the geometric special fiber, and let $D^{b,\rm
Weil}_c(X, \Ql)$ denote the category consisting of objects in the ``derived
category'' $D^b_c(X_{\bar s},\Ql)$ of $\ell$-adic $\Ql$-sheaves, equipped
with a Weil structure.  (See e.g. \cite{GH}, $\S 2$).

We assume that the special fibre $X_s$ of $X$ is stratified by locally
closed subschemes $X_\alpha$, and that each $X_\alpha$ is equidimensional,
of dimension $d_\alpha$, say.

\begin{Def} We say that $\cF \in D^{b,\rm Weil}_c(X, \Ql)$ generically
satisfies the sharp bound with respect to $d$, if for each stratum
$X_\alpha$ there exists an open dense subset $U_\alpha \subseteq X_\alpha$
such that for all $z\in U_\alpha$, the weights on the stalk $\cF_z$ are
$\le d-\dim X_\alpha$.

We say that $\cF$ satisfies the sharp bound with respect to $d$ if for each
$\alpha$ we may take $U_\alpha = X_\alpha$.

We say that $\cF$ (generically) satisfies the sharp bound on its weights,
if it (generically) satisfies the sharp bound w. r. t. $d = \dim X_s$.
\end{Def}

For any $X/S$ such that $X_\eta$ is geometrically integral, we let
$IC(X_\eta)$ denote the (perverse) intersection complex of $X_\eta$, and we
let $^0IC(X_\eta) = IC(X_\eta)(\frac{{\rm dim} X_\eta}{2})$ denote the
Tate-twist of $IC(X_\eta)$ that is a self-dual and weight zero perverse
sheaf on $X_\eta$.  Finally, we define the ``normalized'' nearby cycles
complex in $D^{b,\rm Weil}_c(X,\Ql)$ by $$ ^0R\Psi(IC(X_\eta)) = R\Psi^X(\,
^0IC(X_\eta)).  $$ It is well-known that nearby cycles preserve (middle)
perversity, and so $^0R\Psi(IC(X_\eta))$ is an object of $P_{\rm
Weil}(X,\Ql)$, the full subcategory $D^{b,\rm Weil}_c(X,\Ql)$ whose objects
are perverse. 

Note that it would be more correct to say that the nearby cycles is an
object of the category $D^b_c(X \times_s \eta, \Ql)$, or of its subcategory
$P(X \times_s \eta,\Ql)$ of perverse objects (see e.g. \cite{GH},$\S2,5$).
But as in loc.~cit., by choosing a lift $\Fr_q \in {\rm Gal}(k(\bar
\eta)/k(\eta))$ of the geometric Frobenius generator in ${\rm Gal}(k(\bar
s)/k(s))$, we can and will regard the nearby cycles as being an object in
$P_{\rm Weil}(X,\Ql)$.  This causes no problems in the statement or
application of the proposition below.

\begin{stz}  \label{general_bounds}  Let $X/S$ be an $S$-variety such that
$X_\eta$ is geometrically integral.  Suppose that the special fiber $X_s$
is stratified by locally closed subvarieties $X_\alpha$ which are
equidimensional of dimension $d_\alpha$.  Then $^0R\Psi(IC(X_\eta))$
generically satisfies the sharp bound on its weights.  \end{stz}

%Note that $^0R\Psi(IC(X_\eta))$ is a {\rm mixed} perverse sheaf, as proved
%in \cite{GH}, $\S 10$.  Hence this proposition can be regarded as giving
%information about the weight filtration of $^0R\Psi(IC(X_\eta))$ in the
%category of mixed perverse sheaves.

\noindent {\em Proof.}  By the calculations of Rapoport-Zink \cite{RZ}, the
sharp bound holds in the case where $X/S$ is proper strictly semi-stable,
in the terminology of de Jong \cite{deJ}.  We will reduce the proposition
to that case.  

The key ingredient in the reduction is the following result of de Jong
(\cite{deJ}, Theorem 6.5).

\begin{thm} [de Jong] \label{alterations} Let $X/S$ be an $S$-variety.
Then there exists a trait $S' = (S',s',\eta')$ finite over $S$, an
$S'$-variety $X'$ for which there is an alteration $X' \rightarrow X$ of
$S$-varieties, and an open immersion $j: X' \hookrightarrow \overline{X'}$
of $S'$-varieties such that $\overline{X'}$ is a proper strictly
semi-stable $S'$-variety: $$\xymatrix{ \overline{X'} \ar[d] & X' \ar[l]_j
\ar[d] \ar[r] & X \ar[d] \\ S' & S' \ar[l]_{\rm id} \ar[r] & S.} $$
\end{thm}

We will make use of this theorem in exactly the same manner as \cite{GH},
$\S 10$, where we proved that $^0R\Psi(IC(X_\eta))$ is mixed. 

Let $r: X_{S'} \rightarrow X$ be the projection morphism defined by the
diagram in Theorem \ref{alterations}.  Note that $^0IC(X_{\eta'}) = \,
^0IC(X_\eta)_{\eta'}$.  By invariance of nearby cycles under change of
traits (([SGA 4\,\,1/2, Thm.~finitude 3.7]), it then follows that $r^*
R\Psi^{X}(\, ^0IC(X_\eta)) = R\Psi^{X_{S'}}(\, ^0IC(X_{\eta'}))$, as
objects of $P_{\rm Weil}(X_{S'},\Ql)$.  We have used here that $X_{S'}$ and
$X$ have ``the same'' geometric special fiber $X_{\bar s}$, since $k(s)$ is
finite.  Thus, to bound the weights of $^0R\Psi(IC(X_\eta))$ generically on
each stratum $X_\alpha$, it is enough to do so for $R\Psi^{X_{S'}}(\,
^0IC(X_{\eta'}))$ on each stratum $X_{\alpha,s'}$, the base-change of
$X_\alpha$ via $s' \rightarrow s$.  

We may therefore replace $X/S$ with $X_{S'}/S'$.  In effect, this means we
may assume $S' = S$ and that we have an alteration $\pi: X' \rightarrow X$
of $S$-varieties and an open immersion $j:X' \hookrightarrow \overline{X'}$
of $S$-varieties such that $\overline{X'}/S$ is a proper strictly
semi-stable $S$-variety.

Let us denote simply by $\pi$ the induced morphisms $\pi: X'_{\eta}
\rightarrow X_{\eta}$ and $\pi: X'_{s} \rightarrow X_s$ on generic and
special fibers, and on the corresponding geometric fibers.   For any object
$\mathcal K$ of $D^b_c(X_{\bar \eta},\Ql)$ or $D^b_c(X_{\bar s},\Ql)$, let
$^p{\rm H}^i\mathcal K$ (resp. ${\rm H}^i\mathcal K$) denote the $i$-th
cohomology sheaf of $\mathcal K$ for the perverse (resp. standard)
$t$-structure.  Also, let $^p\tau^{\geq i}$ be the perverse truncation
functor.  Let $^p\pi_* := \, ^p{\rm H}^0\pi_*$ denote the perverse version
of the push-forward functor $\pi_* := R\pi_*$ on derived categories.  Since
$\pi$ is proper, $\pi_* = \pi_!$ and hence $^p\pi_* = \, ^p\pi_!$.   

Let $\Lambda^{X'}_{\eta}$ denote the constant sheaf $\Ql$ on $X'_{\eta}$.
For $n = {\rm dim}(X_\eta)$ let $A^{X'} =
\Lambda^{X'}_{\eta}[n](\frac{n}{2})$.  The smoothness of $X'_{\eta}$ means
that $A^{X'}$ is a self-dual perverse sheaf on $X'_{\eta}$ which is pure of
weight $0$.  

In this circumstance, we know by \cite{GH}, Lemma 10.7, that
$^0IC(X_{\eta})$ is a subquotient of the perverse sheaf $^p\pi_*(A^{X'})$.
By the exactness of the functor $R\Psi^{X}$, the perverse sheaf
$R\Psi^{X}(\, ^0IC(X_{\eta}))$ is a subquotient of $R\Psi^{X}(\,
^p\pi_*A^{X'}) = \, ^p \pi_*R\Psi^{X'}(A^{X'})$ in the category $P_{\rm
Weil}(X,\Ql)$.  Thus, it is enough to prove the sharp bound holds for the
complex $^p\pi_!R\Psi^{X'}(A^{X'})_z$, whenever $z \in U_\alpha$ for some
open dense subset $U_\alpha \subseteq X_\alpha$. 

To bound the weights of $^p\pi_!R\Psi^{X'}(A^{X'})_z$, we claim it is
enough to bound the weights of $\pi_!R\Psi^{X'}(A^{X'})_z$.   To see this,
note that $\pi_!A^{X'}$ is a pure complex of ``geometric origin'' on
$X_{\eta}$, in the terminology of \cite{BBD}, $\S 5, 6$.   Hence for each
$i$ the distinguished triangle $$ \xymatrix{ ^p{\rm H}^i \pi_! A^{X'}
\ar[r] & \, ^p\tau^{\geq i}\pi_! A^{X'} \ar[r] &  \, ^p\tau^{\geq i+1}
\pi_! A^{X'}} $$ becomes a direct sum over $\bar{\eta}$
(loc.~cit.~Thm.~5.4.5, 6.2.5).  Applying the functor $R\Psi^{X'}(\cdot)_z$
we get a distinguished triangle $$ \xymatrix{ ^p{\rm
H}^i\pi_!R\Psi^{X'}(A^{X'})_z \ar[r] & \, ^p\tau^{\geq i}
\pi_!R\Psi^{X'}(A^{X'})_z \ar[r] & \, ^p\tau^{\geq
i+1}\pi_!R\Psi^{X'}(A^{X'})_z} $$ which becomes a direct sum once we forget
the Galois actions.  Thus, for each $j$ we have a short exact sequence of
$\ell$-adic Galois modules $$ \xymatrix{ 0 \ar[r] & {\rm H}^j(\, ^p{\rm
H}^i\mathcal K_z) \ar[r] & \, {\rm H}^j(\, ^p\tau^{\geq i} \mathcal K_z)
\ar[r] & \, {\rm H}^j(\, ^p\tau^{\geq i+1}\mathcal K_z) \ar[r] & 0,} $$
where for brevity we have written $\mathcal K_z$ in place of
$\pi_!R\Psi^{X'}(A^{X'})_z$.  Note that for $i <\!\!< 0$ the middle term is
${\rm H}^j(\pi_!R\Psi^{X'}(A^{X'})_z)$.  It now follows easily by ascending
induction on $i$ that if the sharp bound generically holds for the weights
of $\pi_!R\Psi^{X'}(A^{X'})_z$ then the sharp bound generically holds for
each $^p{\rm H}^i\pi_!R\Psi^{X'}(A^{X'})_z$, and in particular for
$^p\pi_!R\Psi^{X'}(A^{X'})_z$.

Now returning to the diagram in Theorem \ref{alterations}, let
$\Lambda^{\overline{X'}}_{\eta}$ denote the constant sheaf $\Ql$ on
$\overline{X'}_{\eta}$, with corresponding perverse Tate-twist
$A^{\overline{X'}} := \Lambda^{\overline{X'}}_{\eta}[n](\frac{n}{2})$.  Let
$D_i$ be the irreducible components on the special fiber of $X'$.  The
formula $$ j^*R\Psi^{\overline{X'}}(A^{\overline{X'}}) = R\Psi^{X'}(A^{X'})
$$ and the results of Rapoport-Zink \cite{RZ} for the proper strictly
semi-stable model $\overline{X'} \rightarrow S$ show that the weights of
$R\Psi^{X'}(A^{X'})$ on any intersection $\cap_{i \in J}D_i$ are $\leq {\rm
dim} (X'_{\bar {s}}) - {\rm dim}(\cap_{i \in J}D_i)$.  Here $J$ denotes any
finite subset of the set indexing the components $D_i$.

Now we can complete the proof of the proposition.  Since $\pi : X'_{s}
\rightarrow X_{s}$ is proper, the union of the sets $X_{\alpha} \cap
\pi(\cap_{i \in I}D_i)$, where $I$ ranges over all subsets of indices such
that ${\rm dim}(\cap_{i \in I}D_i) < d_\alpha$, is a proper closed subset
of $X_{\alpha}$.  Moreover, since $X_{\alpha}$ is equidimensional the
complement $U_{\alpha} \subseteq X_{\alpha}$ of this closed subset is open
and {\em dense}. For $z \in U_{\alpha}$ the fiber $\pi^{-1}(z)$ does not
meet any intersection of irreducible components in $X'_{s}$ of form
$\cap_{i \in I} D_i$ whose dimension is $< d_\alpha$.  It follows that the
weights of $R\Psi^{X'}(A^{X'})$ on the fiber $\pi^{-1}(z)$ are $\leq {\rm
dim}(X'_{\bar{s}}) - d_\alpha$, by the above remarks.  Hence for $z \in
U_{\alpha}$ the weights of $\pi_!R\Psi^{X'}(A^{X'})_z =
R\Gamma_c(\pi^{-1}(z), R\Psi^{X'}(A^{X'}))$ are $\leq {\rm
dim}(X'_{\bar{s}}) - d_\alpha$, and we are done.  \qed

\subsection{Application to the affine flag manifold}

Of course, the generic sharp bound on the weights is most interesting in
cases where the weights are actually constant along the strata, in
particular when one is dealing with sheaves equivariant under some group
action, where the stratification is given by the orbits.  This is the
situation we considered in \cite{GH}: in either the $p$-adic or the
function-field setting, the sheaves $R\Psi_\mu$ we consider are
$\cB$-equivariant hence constant along the $\cB$-orbits $\cB_w$, for $w \in
\widetilde{W}$; see loc.~cit.~$\S 2$ and \cite{HN1},\cite{Ga}.  Hence we
have the following corollary of Proposition \ref{general_bounds}, which
immediately proves Theorems \ref{nearby_bounds} and \ref{nearby.poly}.

\begin{kor} Since $R\Psi_\mu(\frac{\ell(\mu)}{2})$ is $\cB$-equivariant, it
satisfies the sharp bound on its weights with respect to $d = \ell(\mu) =
{\rm dim}({\mathcal Q}_\mu)$ and the stratification by $\cB$-orbits.
\end{kor}

\section[Coefficients in Wakimoto functions]{Bounding degrees of
coefficients in Wakimoto functions} \label{wakimoto_section}

In this section we give a combinatorial result which proves Theorem
\ref{wakimoto.poly} and which also yields ---via the Kottwitz conjecture
(\ref{eq:Kottwitz_conj})--- a different proof of Theorem \ref{nearby.poly}.

By Lemma \ref{weights_v_degrees}, Theorem \ref{wakimoto.poly} can be
reformulated as follows.  

\begin{thm} \label{wakimoto.poly.reformulated} For each $w \leq uv$, the
function $\Tr(\Fr_q, (\widetilde{M}_{u,v})_w)$ is a polynomial in $q$
having degree $\leq \ell(uv) - \ell(w)$.  \end{thm}

We will need some further notation.  Recall we consider a split connected
reductive group $G$ over a field $k$.  Let us fix a choice of Borel
subgroup $B \subset G$ and maximal torus $T \subset B$.  Let $\cB \subset
G(k\xt)$ be the Iwahori subgroup whose ``reduction modulo $t$'' is $B$.
The data $(G,B,T)$ determine a (based) root system and a set of simple
reflections $S$ generating the Coxeter group $(W,S)$, where $W = N_GT/T$ is
the finite Weyl group.  Let $w_0$ denote the unique longest element in $W$.
Let us choose an origin in the apartment corresponding to $T$ and a Weyl
chamber ${\rm Ch}(B)$ corresponding to $B$.  Let $A$ denote the {\em base
alcove}, that is, the unique alcove in the given apartment whose closure
contains the origin and which lies in ${\rm Ch}(B)$.  Let $\bar{B}$ be the
unique {\em opposite} Borel subgroup, that is, the one such that $B \cap
\bar{B} = T$.    

Let $W_{\rm aff}$ denote the affine Weyl group of $G$.  It is a Coxeter
group with generators $S_{\rm aff}$, the set of simple reflections through
the walls of the {\em opposite} alcove $w_0(A)$ of the base alcove $A$: if
$G$ is almost simple with $B$-positive simple roots $\alpha_i$ ($1 \leq i
\leq l$) and $B$-highest root $\widetilde{\alpha}$, then $w(A_0)$ is the
alcove $$ w_0(A) = \{ x \in X_*(T) \otimes {\mathbb R} ~ | ~ \mbox{$0 <
\langle -\alpha_i,x \rangle < 1$ for all $i$, and $\langle
-\widetilde{\alpha}, x \rangle < 1$} \}.  $$ In this case the simple affine
reflections are the reflections $s_{-\alpha_i} = s_{\alpha_i}$ and $s_0 :=
s_{\widetilde{\alpha} + 1} =
t_{-\widetilde{\alpha}^\vee}s_{\widetilde{\alpha}}$. (See also \cite{GH},
section 2.2.)

The alcove fixed by $\cB$ is the {\em opposite} alcove $w_0(A)$.  This
results from our convention that we embed $\widetilde{W}$ into $G(k\xT)$
using the identification of $\nu \in X_*(T)$ with $\nu(t) \in G(k\xT)$. 
  
The $\bar{B}$-positive affine roots $\alpha + k$ are those taking positive
values on $w_0(A)$.  Equivalently, $\alpha + k$ is $\bar{B}$-positive if
and only if either $k \geq 1$, or $k = 0$ and $\alpha$ is
$\bar{B}$-positive.

The set of affine roots carries a left action by the group $\widetilde{W}$
by the rule $$ x \cdot (\alpha + k) = (\alpha + k ) \circ x^{-1}.  $$

\begin{lem} \label{relating_orders} Let $x \in \widetilde{W}$ and suppose
$\alpha + k$ is a $\bar{B}$-positive affine root, with corresponding
reflection $s = s_{\alpha + k}$.  Let $\leq$ denote the Bruhat order on
$\widetilde{W}$ defined using $S_{\rm aff}$.  Then \begin{enumerate} \item
[(a)] $x < xs$ if and only if $x \cdot (\alpha + k)$ is $\bar{B}$-positive;
\item [(b)]$ x < sx$ if and only if $x^{-1} \cdot (\alpha + k)$ is
$\bar{B}$-positive.  \end{enumerate} \end{lem}

\noindent {\em Proof.}  It is enough to prove (b).  Write $A^- := w_0(A)$.
We note that $x < sx$ if and only if $xA^-$ and $A^-$ lie on the same side
of the hyperplane $\alpha + k = 0$, which happens if and only if $(\alpha +
k) \circ x$ takes positive values on $A^-$, i.e., $x^{-1} \cdot (\alpha +
k)$ is  $\bar{B}$-positive.  \qed

We will also need one more lemma.

\begin{lem} \label{sublemma} Suppose $t = s_\beta$, where $\beta$ is a
$\bar{B}$-positive affine root.  Let $u,v \in \widetilde{W}$ and suppose $u
< ut$.  Then $$ uv < utv \Longleftrightarrow v < tv.  $$ \end{lem}

\noindent{\em Proof.}  Note first that $utu^{-1} = s_{u \cdot \beta}$ and
$u \cdot \beta$ is $\bar{B}$-positive, by Lemma \ref{relating_orders}.
Then we have \begin{align*} uv < utv &\Longleftrightarrow uv < s_{u\cdot
\beta} uv \\ &\Longleftrightarrow \mbox{$(uv)^{-1} \cdot (u \cdot \beta)$
is $\bar{B}$-positive} \\ &\Longleftrightarrow \mbox{$v^{-1} \cdot \beta$
is $\bar{B}$-positive} \\ &\Longleftrightarrow v < tv, \end{align*} again
by Lemma \ref{relating_orders}.  \qed

\medskip

For elements $u, v \in \widetilde{W}$, we consider the {\em Wakimoto}
functions $\widetilde{T}_u \widetilde{T}_{v^{-1}}^{-1}$, see \cite{GH},
section 7.  Here for $x \in \widetilde{W}$, the symbol $T_x$ denotes the
standard generator $\rm char(\cB x \cB)$ in the Iwahori-Hecke algbra
$C_c(\cB\backslash G(\F_q\xT)/\cB)$ and $\widetilde{T}_x := q_x^{-1/2}T_x$.
Our goal is to prove the following improved version of
loc.~cit.~Proposition 7.6 (where the same conclusion is obtained under the
additional hypothesis that $\widetilde{T}_u \widetilde{T}^{-1}_{v^{-1}}$
has a minimal expression).  The indeterminate $Q$ that appears here is
defined by the relation $q^{1/2}Q = 1-q$.

\begin{stz} \label{Wakimoto_fcns_bound} If $\widetilde{T}_u
\widetilde{T}^{-1}_{v^{-1}} = \sum_x R^u_{x,v}(Q) \, \widetilde{T}_x$, then
$R^u_{x,v}(Q)$ is a polynomial in $Q$ with integer coefficients, and $$
{\rm deg}_Q R^u_{x,v}(Q) \leq \ell(uv) - \ell(x).  $$ \end{stz}

\noindent{\em Proof.}  We may assume $u,v \in W_{\rm aff}$.  Consider
reduced expressions $u = t_1 \cdots t_k$ and $v = s_1 \cdots s_r$, where
$t_j, s_i$ belong to $S_{\rm aff}$.  Write $\tilde{t}_j$ resp.
$\tilde{s}_i^{-1}$ as short-hand for $\widetilde{T}_{t_j}$ resp.
$\widetilde{T}^{-1}_{s_i}$.  We will prove the desired bound $$ {\rm deg}_Q
\, \tilde{t}_1 \cdots \tilde{t}_k \cdot \tilde{s}_1^{-1} \cdots
\tilde{s}_r^{-1} (x) \leq \ell(uv) - \ell(x), $$ by induction on $k$, the
case $k =0$ being obvious.  We assume the bound holds for $k \geq 0$ and
deduce it for $k+1$.

Let $t = s_\beta$, where $\beta$ is a $\bar{B}$-positive simple affine
root, and suppose $u < ut = t_1 \cdots t_kt$.  

\medskip

\noindent {\em Case 1:}  We have $v < tv$  (so by Lemma \ref{sublemma}, $uv
< utv$).  Then because $\tilde{t}^{-1} = \tilde{t} + Q$, we have $$
\tilde{t}_1 \cdots \tilde{t}_k\tilde{t} \tilde{s}_1^{-1} \cdots
\tilde{s}_r^{-1} (x) = \tilde{t}_1 \cdots \tilde{t}_k \tilde{t}^{-1}
\tilde{s}_1^{-1} \cdots \tilde{s}_r^{-1}(x) - Q \tilde{t}_1 \cdots
\tilde{t}_k \tilde{s}_1^{-1} \cdots \tilde{s}_r^{-1}(x).  $$ Now by
induction the ${\rm deg}_Q$ of the first term on the right is bounded by
$\ell(utv) - \ell(x)$.  Similarly the ${\rm deg}_Q$ of the second term is
bounded by $1 + \ell(uv) - \ell(x) \leq \ell(utv) - \ell(x)$.

\medskip

\noindent {\em Case 2:} We have $tv < v$ (so by Lemma \ref{sublemma}, $utv
< uv$).  Then we can write a reduced expression for $v$ in the form $v =
ts_2 \cdots s_r$, and thus we are considering $$ \tilde{t}_1 \cdots
\tilde{t}_k \tilde{t} \tilde{t}^{-1} \tilde{s}_2^{-1} \cdots
\tilde{s}_r^{-1} (x) = \tilde{t}_1 \cdots \tilde{t}_k \tilde{s}_2^{-1}
\cdots \tilde{s}_r^{-1} (x)$$ and clearly our induction hypothesis implies
that the ${\rm deg}_Q$ is bounded by $\ell(utv) - \ell(x)$, as desired.  

The fact that $R^u_{x,v}$ is a {\em polynomial} in $Q$ with integer
coefficients can also be proved by induction on $k = \ell(u)$ in a similiar
way, or by using the Iwahori-Matsumoto relations in the Hecke algebra (see
\cite{H01}).  \qed

\bigskip

Let us now state explicitly what Proposition \ref{Wakimoto_fcns_bound} says
about the weights of the Wakimoto sheaves $\widetilde{M}_{u,v}$. Recall
that for $u, v \in \widetilde{W}$, the Wakimoto sheaf is defined as \[
M_{u,v} = j_{u!} \Ql[\ell(u)] \star j_{v*}\Ql[\ell(v)], \] where $j_u$,
$j_v$ are the inclusions of the Schubert cells corresponding to $u$, $v$
into the affine flag variety, where $j_{u!}$, $j_{v*}$ denote the derived
functors, and where $\star$ is the convolution product of perverse sheaves.
The function associated to $M_{u,v}$ via the sheaf-function dictionary is
\[ \varepsilon_u \varepsilon_v q_v T_u T_{v^{-1}}^{-1} = \varepsilon_u
\varepsilon_v q_u^{1/2} q_v^{1/2}\sum_x q_x^{-1/2} R^u_{x,v} T_x, \] and so
the function associated to $\widetilde{M}_{u,v}$ is $\varepsilon_u
\varepsilon_v q_{uv}^{1/2}\sum_x q_x^{-1/2} R^u_{x,v} T_x$.

From \cite{GH}, $\S 7$ we know that for any power $Q^i$ which appears with
non-zero coefficient in $R^u_{x,v}(Q)$, the integer $i$ has the same parity
as $\ell(uv) - \ell(x)$.  This together with the identity $q^{1/2}Q=1-q$
and the bound ${\rm deg}_Q R^u_{x,v}(Q) \leq \ell(uv) - \ell(x)$ of
Proposition \ref{Wakimoto_fcns_bound} gives the following result, which is
a reformulation of Theorems \ref{wakimoto.poly.reformulated} and
\ref{wakimoto.poly}.

\begin{kor} \label{wakimoto_bounds} Let $u, v, x \in \widetilde{W}$. Then
$\Tr(\Frob_q, (\widetilde{M}_{u,v})_x)$ is a polynomial in $q$ of degree
$\le \ell(uv) - \ell(x)$. The sheaf
$\widetilde{M}_{u,v}(\frac{\ell(uv)}{2})$ satisfies the sharp bound with
respect to $d=\ell(uv)$ and the stratification by $\cB$-orbits.  \end{kor}

We remark that via the Kottwitz conjecture stated below, this corollary can
be used to give another proof of Theorem \ref{nearby.poly}.  The Kottwitz
conjecture is the following identity in the Iwahori-Hecke algebra $C_c(\cB
\backslash G(\F_q\xT)/\cB)$.  In the function-field setting it was proved
by Gaitsgory \cite{Ga}, whose ideas were adapted to prove it in the
$p$-adic setting in \cite{HN1}.

\begin{equation}\label{eq:Kottwitz_conj} \Tr(\Fr_q, R\Psi_\mu) =
\varepsilon_\mu q_\mu^{1/2} \sum_{\lambda \in \Omega(\mu)} m_\mu(\lambda)
\Theta_\lambda.  \end{equation}

The original form concerned the {\em semi-simple} trace of Frobenius, but
by virtue of the unipotence of the inertia action on $R\Psi_\mu$ which was
proved by Gaitsgory \cite{Ga} (see also \cite{GH} $\S 5$), we may replace
the semi-simple trace $\Tr^{ss}$ with the usual trace, once we have chosen
a lift $\Fr_q$ of geometric Frobenius, as we have already done.  See the
discussion of this point in \cite{GH}, $\S 5$.

Let us explain the notation on the right hand side.  We let $\Omega(\mu)$
denote the set of coweights $\lambda \in X_*(T)$ such that for every $w \in
W$, $\mu - w\lambda$ is a sum of $B$-positive coroots.  For such a coweight
$\lambda$, write it as $\lambda = \lambda_1 - \lambda_2$, where $\lambda_i$
is $B$-dominant for $i= 1,2$.  Then following Bernstein we define $$
\Theta_\lambda = \widetilde{T}_{t_{\lambda_1}}
\widetilde{T}^{-1}_{t_{\lambda_2}}.  $$ Let $G^\vee$ be the dual group for
$G$.  The number $m_\mu(\lambda)$ is the multiplicity with which the weight
$\lambda$ appears in the character $E_\mu$ of the simple $G^\vee$-module
with highest weight $\mu$.

Note that $\ell(\lambda) \leq \ell(\mu)$ and the difference $\ell(\mu) -
\ell(\lambda)$ is an even integer.  It is now clear that Theorem
\ref{nearby.poly} can be deduced, via Lemma \ref{weights_v_degrees}, by
applying Corollary \ref{wakimoto_bounds} with $u = t_{\lambda_1}$ and $v =
t_{-\lambda_2}$ for all $\lambda \in \Omega(\mu)$ and by using
(\ref{eq:Kottwitz_conj}).

\bigskip \bigskip

Ulrich G\"{o}rtz: Mathematisches Institut der Universit\"{a}t Bonn,
Beringstr. 6, 53115 Bonn, Germany.  Email: ugoertz@math.uni-bonn.de

\bigskip

Thomas J. Haines: Mathematics Department, University of Maryland, College
Park, MD 20742-4015, USA.  Email: tjh@math.umd.edu

\end{document}